\documentclass{amsproc}
\usepackage{amsmath,amssymb,amsfonts}

\newtheorem{theorem}{Theorem}[section]
\newtheorem{lemma}[theorem]{Lemma}
\newtheorem*{corollary}{\textnormal{Corollary}}

\theoremstyle{definition}
\newtheorem{definition}[theorem]{Definition}
\newtheorem{example}[theorem]{Example}

\theoremstyle{remark}
\newtheorem{remark}[theorem]{Remark}
\newcommand{\alg}{\operatorname{Alg}}
\newcommand{\lat}{\operatorname{Lat}}

\numberwithin{equation}{section}

\begin{document}

\title{Rank Preserving Maps on CSL Algebras}

\author{Jaedeok Kim}
\address{Department of Mathematical Computing and Information Sciences, Jacksonville State University, Jacksonville, AL 36265}
\email{jkim@jsu.edu}

\author{Robert L. Moore}
\address{Department of Mathematics, University of Alabama,
Tuscaloosa, AL 35487 }\email{rmoore@gp.as.ua.edu}

\date{January 1, 1994 and, in revised form, June 22, 1994.}



\begin{abstract}
We give a description of a weakly continuous rank preserving map
on a reflexive algebra on complex Hilbert space with commutative
completely distributive subspace lattice. We show that the
implementation of a rank preserving map can be described by the
combination of two different types of maps. We also show that a
rank preserving map can be implemented by only one type if the
corresponding lattice is irreducible. We present some examples of
both types of rank preserving map.
\end{abstract}

\maketitle



\section{PRELIMINARIES}
It has been of special interest to study the linear maps between
two
 nonself\-adjoint operator algebras over the last several decades.
 These maps include isomorphisms, isometries and rank preserving
 maps etc. In \cite{MR33:4703}, Ringrose proved that an isomorphism between two
 nest algebras must be spatially implemented by making effective
 use of rank one operators. In \cite{MR87k:47103}, Gilfeather and Moore showed
 that an isomorphism acting between two CSL algebras with
 completely \mbox{distributive} lattices need not be spatially or
 quasi-spatially implemented. But they also showed that if an
 isomorphism preserves the rank of all finite rank operators, then
 it must be quasi-spatially implemented. In \cite{MR17:09095}, Panaia has shown that
 if $\mathfrak{A}_i=\alg{\mathfrak{L}}_i$ where the $\mathfrak{L}_i$ are finite distributive
 subspaces lattices, then every rank preserving algebraic isomorphism of $\mathfrak{A}_1$ onto
 $\mathfrak{A}_2$ is quasi-spatial. The characterization of
 the rank preserving maps on nest algebras was done by the Chinese
 mathematicians, Shu-Yun Wei and Sheng-Zhao Hou \cite{MR2000m:47099}. Part of their results showed that
some rank preserving maps $\Phi$ on a nest algebra can be
described as $\Phi(T)=ATB^*$ for some $A,B \in B(\mathcal{H})$.
Use of rank one operators was an integral part of their proof. In
this paper, we shall describe the rank preserving maps on
completely distributive commutative subspace lattice algebras
(CDCSL algebras).

In this paper we show that a rank preserving map $\Phi$ on a
completely distributive algebra, CDC algebra, with certain
condition can be described as $\Phi(A)=UAV^*$ for some densely
defined linear transformations $U$ and $V$.

Let $\mathcal{H}$ be a complex separable Hilbert space. A
\emph{subspace lattice} $\mathfrak{L}$ is a strongly closed
lattice of orthogonal projections on $\mathcal{H}$, containing $0$
and $I$. If $\mathfrak{L}$ is a subspace lattice,
$\alg{\mathfrak{L}}$ denotes the algebra of all bounded operators
on $\mathcal{H}$ that leave invariant every projection in
$\mathfrak{L}$. $\alg{\mathfrak{L}}$ is a weakly closed subalgebra
of $B(\mathcal{H})$, the algebra of all bounded operators on
$\mathcal{H}$. Dually, if $\mathfrak{A}$ is a subalgebra of
$B(\mathcal{H})$, then $\lat{\mathfrak{A}}$ is the lattice of all
projections invariant for each operator in $\mathfrak{A}$. An
algebra $\mathfrak{A}$ is reflexive if
$\mathfrak{A}=\alg\lat{\mathfrak{A}}$ and a lattice $\mathfrak{L}$
is reflexive if $\mathfrak{L}=\lat\alg{\mathfrak{L}}$. A lattice
is a \emph{commutative subspace lattice}, or CSL, if each pair of
projections in $\mathfrak{L}$ commute; $\alg{\mathfrak{L}}$ is
then called a CSL algebra. All lattices in this paper will be
commutative.

In \cite{MR51:1420}, Arveson showed that every commutative
subspace lattice is reflexive. A totally ordered (and hence
commutative) subspace lattice is a nest and the associated algebra
is a nest algebra. We use the convention that a subspace is
identified with the orthogonal projection onto the subspace. Thus
$E \subseteq F$ is the same as $E \leq F$.

We have two operations between subspaces: meet($\wedge$) and
join($\vee$).

Given any family $\{E_{\alpha}\}_{\alpha \in I}$ of subspaces of a
Hilbert space $\mathcal{H}$, $\wedge_{\alpha\in I}E_{\alpha}$
denotes the greatest subspace contained in each $E_{\alpha}$ and
$\vee_{\alpha\in I}E_{\alpha}$ denotes the smallest subspace
containing each $E_{\alpha}$.

We now restrict our attention to a special type of lattice:
completely distributive lattice. One of the advantages of
completely distributive lattices over other lattices is the
abundance of the rank one operators in the associated algebra;
these operators will be building blocks in the subsequent
discussion. A \emph{completely distributive lattice} is a complete
lattice which satisfies distributive laws expressed for families
of arbitrary cardinality. To be more precise, let $\mathfrak{L}$
be a complete lattice and $I$ be an arbitrary index set. For each
$\alpha \in I$, let $J_{\alpha}$ also be an arbitrary index set,
and for each $\beta \in J_{\alpha}$, let $E_{\beta \alpha}$ denote
an element of $\mathfrak{L}$. $\Pi$ will denote the Cartesian
product of all the $J_{\alpha}$, i.e. the collection of all choice
functions $\varphi:I \rightarrow\cup_{\alpha \in I}J_{\alpha}$
satisfying $\varphi(\alpha) \in J_{\alpha}$, for all $\alpha$.
$\mathfrak{L}$ is completely distributive if it satisfies the
following two identities.

\begin{enumerate}
\item $\wedge_{\alpha \in I}(\vee_{\beta \in
J_{\alpha}}E_{\beta\alpha})=\vee_{\varphi \in \Pi}(\wedge_{\alpha
\in I}E_{\varphi(\alpha)\alpha})$

\item $\vee_{\alpha \in I}(\wedge_{\beta \in
J_{\alpha}}E_{\beta\alpha})=\wedge_{\varphi \in \Pi}(\vee_{\alpha
\in I}E_{\varphi(\alpha)\alpha})$
\end{enumerate}

If $\mathfrak{L}$ is completely distributive and commutative, we
will call $\alg{\mathfrak{L}}$ a CDC algebra. For further
discussion about the characteristics of CDC algebra, we need to
introduce the notion of rank one operators. We will let $x\otimes
y^{\ast}$ denote the rank one operator defined on $\mathcal{H}$ by
$(x\otimes y^{\ast})(f)=\langle f,y \rangle x$.

 There are plenty of rank one operators and finite rank
operators in a nest algebra. Moreover, these operators have played
a central role in the theory of nest algebra. However, there are
many examples of commutative subspace lattices whose corresponding
algebras do not contain any rank one operators \cite{MR84h:47048}.
But as we discussed earlier, the complete distributivity of a
lattice guarantees a reasonable supply of rank one operators.
Furthermore, one of the surprising characteristics of the CDC
algebra is stated in the following lemma. This lemma is due to
Laurie and Longstaff \cite{MR85b:47052}.

\begin{lemma}\label{wk}
Let $\mathfrak{L}$ be a commutative completely distributive
subspace lattice on a separable Hilbert space $\mathcal{H}$. Let
$R_{\mathfrak{L}}$ denote the linear span of the rank one
operators in $\alg \mathfrak{L}$. Then $R_{\mathfrak{L}}$ is dense
in $\alg\mathfrak{L}$ in any of the weak, strong, ultraweak or
ultrastrong topologies.
\end{lemma}

It is an interesting question to ask when $x \otimes y^{\ast}$
belongs to $\alg\mathfrak{L}$. The following lemma, due to
Longstaff \cite{MR53:1294}, is the answer to this question, and
the use of this lemma is essential.

\begin{lemma}\label{lsl1}
The operator $x \otimes f^{\ast} $ belongs to $\alg \mathfrak{L}$
if and only if there is a projection $E \in \mathfrak{L}$ such
that $x \in E $ and $f \in E_{-}^{\perp}$.
\end{lemma}

Let $\mathfrak{L}$ be a subspace lattice. For $M \in
\mathfrak{L}$, define $M_{-}$ by $M_{-} = \vee \{N|M \neq N, N \in
\mathfrak{L}\}$. We end this section with a couple of lemmas that
will be used repeatedly.

\begin{lemma}\label{cdl1}
Let $\mathfrak{L}$ be a commutative completely distributive
lattice. Then
\[\vee \{N | N \in \mathfrak{L}, N_{-} \neq I \}=I \]
 and
\[\vee \{N_{-}^{\perp}| N \in \mathfrak{L}, N \neq 0  \}=I \]
\end{lemma}

\begin{lemma}\label{ul}
Let $\mathfrak{L}$ be a subspace lattice and $\{E_k\}_{k=1}^n$ be
a finite subset of $\mathfrak{L}$. Then $(\vee_{k=1}^n
E_k)_-=\vee_{k=1}^n (E_k)_- $ and $(\wedge_{k=1}^n E_k)_- \leq
\wedge_{k=1}^n (E_k)_- $
\end{lemma}
\begin{proof}
It is sufficient to prove that the statements are true in case
$n=2$. Let $E$ and $F$ denote two elements in $\mathfrak{L}$.
Since $E \leq E \vee F$ and $F \leq E \vee F$, we have $E_- \leq
(E \vee F)_-$ and $F_- \leq (E \vee F)_-$. These imply that \[ E_-
\vee F_- \leq (E \vee F)_-. \]

 On the other hand, if $K \in \{G
\in \mathfrak{L} \, | \, G \ngeq E \vee F \}$, then $K \in \{G \in
\mathfrak{L} \, | \, G \ngeq E\}$ or $K \in \{G \in \mathfrak{L}
\, | \, G \ngeq F\}$, so $K \leq E_-$ or $K \leq F_-$. These imply
that $K \leq E_- \vee F_-$. Hence $ (E \vee F)_- \leq E_- \vee
F_-.$ This proves $(1)$.

 Let $E$ and $F$ be two elements in
$\mathfrak{L}$. Since $E \wedge F \leq E$ and $E \wedge F \leq F$,
$(E \wedge F)_- \leq E_-$ and $(E \wedge F)_- \leq F_-$. Therefore
$(E \wedge F)_- \leq E_- \wedge F_-$. The results follow.
\end{proof}

Now let's discuss two important types of linear maps on $\alg
\mathfrak{L}$.

\begin{definition}
Let $\mathfrak{L}$ be a commutative subspace lattice.
\begin{enumerate}
\item A linear map $\Phi:\alg\mathfrak{L}\rightarrow
\alg\mathfrak{L}$ is called an \emph{isomorphism} if $\Phi$ is a
bijection and a multiplication preserving map. \item A linear map
$\Phi:\alg \mathfrak{L}\rightarrow \alg\mathfrak{L}$ is called a
\emph{rank preserving map} if $rank(\Phi(A))=rank(A)$ for each
finite rank operator $A$ in $\alg\mathfrak{L}$.
\end{enumerate}
\end{definition}

The automatic norm continuity of an isomorphism is proved by
Gilfeather and Moore in \cite{MR87k:47103}. In order to discuss
more results about isomorphisms, we need the following
terminology. An isomorphism $\Phi$ is said to be \emph{spatially
implemented} if there is a bounded invertible operator $T$ on
$\mathcal{H}$ so that $\Phi(A)=TAT^{-1}$ for all $A \in
\alg(\mathfrak{L})$. An isomorphism is said to be
\emph{quasi}-\emph{spatially implemented} if there is a one-to-one
operator with dense domain $\mathcal{D}$ so that $\Phi(A)Tf=TAf$
for all $A \in \alg\mathcal{L}$ and for all $f \in \mathcal{D}$.

The following theorem characterizes isomorphisms of CSL algebras.
The theorem is due to Gilfeather and Moore \cite{MR87k:47103}.

\begin{theorem}
Let $\mathfrak{L}_1$ and $\mathfrak{L}_2$ be commutative subspace
lattices on Hilbert spaces $\mathcal{H}_1$ and $\mathcal{H}_2$,
respectively, and let $\mathfrak{L}_1$ be completely distributive.
Let $\rho:\alg\mathfrak{L_1}\rightarrow \alg \mathfrak{L_2}$ be an
algebraic isomorphism. The followings are equivalent.
\begin{enumerate}
\item $\rho$ is quasi-spatially implemented by a closed, injective
linear transformation $T:H_1 \rightarrow H_2$ whose range and
domain are dense. \item $\rho$ is a rank-preserving map.
\end{enumerate}
\end{theorem}

The weak continuity of a rank preserving map $\Phi$ will be used
in the proof of the main Theorem. However, we set the assumption
that \emph{each rank preserving map $\Phi$ is weakly continuous}
from this early stage to avoid any misunderstanding.

\section{RANK PRESERVING LINEAR MAPS}
We begin with a lemma which will play a key role in our result.
For notational convenience, we will let $a \sim b$ represent the
linear dependency of two nonzero vectors $a,b$ in $\mathcal{H}$.
It can be easily verified that $a \sim b$ is an equivalence
relation.

Suppose $\mathfrak{L}$ is a commutative subspace lattice and
$\Phi$ is a rank preserving map from $\alg \mathfrak{L}$ to $\alg
\mathfrak{L}$. Let $N$ be an element in $\mathfrak{L}$ with $\dim
N \geq 2, \, \dim N_-^{\perp} \geq 2$. Let $x,y$ be any two
linearly independent vectors in $N$ and $f,g$ be any two linearly
independent vectors in $N_-^{\perp}$. Then we can form four
different rank one operators in $\alg\mathfrak{L}$ as follows: $x
\otimes f^{\ast},y \otimes f^{\ast},x \otimes g^{\ast},y \otimes
g^{\ast}$.

Since $\Phi$ is a rank preserving map on $\alg\mathfrak{L}$, we
can consider another four rank one operators in $\alg
\mathfrak{L}$ which are images of each of the four rank one
operators and write them as follows.
\begin{align*}
  \Phi(x \otimes f^{\ast})&= u \otimes v^{\ast} &&\cdots \cdots (A)    \\
  \Phi(y \otimes f^{\ast})&= p \otimes q^{\ast} &&\cdots \cdots (B)    \\
  \Phi(x \otimes g^{\ast})&= w \otimes z^{\ast} &&\cdots \cdots (C)    \\
  \Phi(y \otimes g^{\ast})&= r \otimes s^{\ast} &&\cdots \cdots (D)
\end{align*}

\begin{lemma} (Four Vectors Lemma)\label{fvl}
 Let $\mathfrak{L}$ be a
commutative subspace lattice and $\Phi$ be a rank preserving map
on $\alg\mathfrak{L}$. Let $N$ be a subspace in $\mathfrak{L}$
with $ \dim N \geq 2, \, \dim N_-^{\perp} \geq 2$. If we consider
the eight rank one operators described above, then either $(1)$ or
$(2)$ holds for all $x,y \in N$ and $f,g \in N_-^{\perp}$.
 \begin{enumerate}
  \item $u \sim w, v \sim q, p \sim r \, and \, z \sim s$
  \item $v \sim z, u \sim p, w \sim r \, and \, q \sim r$
\end{enumerate}
\end{lemma}

\begin{proof}
Let's consider the following six equalities which are obtained by
adding two of $ (A), (B),(C) \, and \, (D) $.

$(A)+(C):\Phi(x \otimes(f+g)^*)=u\otimes v^*+w\otimes z^*$

$(A)+(B):\Phi((x + y)\otimes f^*)=u\otimes v^*+p\otimes q^*$

$(B)+(D):\Phi(y \otimes(f+g)^*)=p\otimes q^*+r\otimes s^*$

$(C)+(D):\Phi((x+y) \otimes g^*)=w\otimes z^*+r\otimes s^*$

$(A)+(D):\Phi(x \otimes f^*+y \otimes g^*)=u\otimes v^*+r\otimes
s^*$

$(B)+(C):\Phi(x \otimes g^*+y\otimes f^*)=p\otimes q^*+w\otimes
z^*$

Since $x\otimes (f+g)^*$ is a rank one operator and $\Phi$ is a
rank preserving map, $u\otimes v^*+w\otimes z^*$ is also a rank
one operator. This fact implies either $u\sim w$ or $v\sim z$.
Similarly we can argue that either $u\sim p$ or $v\sim q$ is true
by the equality $(A)+(B)$, either $p\sim r$ or $q\sim s$ is true
by the equality $(B)+(D)$ and $w\sim r$ or $z\sim s$ is true by
the fourth equality $(C)+(D)$. Apart from the first four
equalities, the map $\Phi$ maps a rank two operator into a rank
two operator in the equality $(A)+(D)$ and$(B)+(C)$, which gives
the following results
\[u \nsim r, v\nsim s, p\nsim w \, \textrm{and} \, q\nsim z .\]
\underline{Case 1:} Now let's suppose $u\sim w$ is true. Then it
is easily deduced that $v\sim q$. Otherwise $u\sim p$ is true
implying $w \sim p$, which contradicts the relation $w \nsim p$.
Similarly, we can show that $p\sim r, z\sim s$ are true, so the
result follows.
\\
\underline{Case 2:} Suppose $v\sim z$. The proof of this case is
almost identical to the proof of Case 1.
\end{proof}

\begin{remark}
If we assume that \emph{Case 1} holds, then $u \nsim p$.
\end{remark}

If we consider the case $\dim N=1$ and $\dim N_{-}^{\perp}\geq 2$,
we have the following two equalities.
\begin{align*}
 \Phi(x \otimes f^{\ast})&= u \otimes v^{\ast} &&\cdots \cdots (E)
    \\
 \Phi(x \otimes g^{\ast})&= p \otimes q^{\ast} &&\cdots \cdots (F)
\end{align*}
where $x \in N$ and $f,g$ are two linearly independent vectors in
$N_{-}^{\perp}$.
 The following lemma is analogous to Lemma
\ref{fvl} for this case.

\begin{lemma}\label{tvl}
 Let $\mathfrak{L}$ be a
commutative subspace lattice and $\Phi$ be a rank preserving map
on $\alg\mathfrak{L}$. Let $N$ be a subspace in $\alg\mathfrak{L}$
with $\dim N=1$ and $\dim N_{-}^{\perp}\geq 2$. If we consider the
two equalities described above, then either $(1)$ or $(2)$ holds
for all $x \in N$ and $f,g \in N_-^{\perp}$.
 \begin{enumerate}
  \item $ u \sim p , \, v \nsim q $
  \item $ u \nsim p , \, v \sim q $
  \end{enumerate}
\end{lemma}

\begin{proof}
By adding the two equalities, the following equality is obtained.

 $(E)+(F): \Phi(x \otimes (f+g)^*)=u \otimes v^* + p \otimes
 q^*$

 From the fact that $\Phi$ maps a rank one operator into a rank one operator,
 it can be easily deduced that either $ u \sim p$ or $ v \sim q$
 is true. Suppose $u \sim p$ and $v \sim q$. Then there exist
 $\lambda, \gamma \in \mathbb{C}$ such that $\lambda u \otimes v^* + \gamma p
 \otimes q^*=0$. Therefore $ \Phi(x \otimes
 (\overline{\lambda} f + \overline{\gamma} g)^*) = \lambda u \otimes v^* + \gamma
 p\otimes q^*=0$ which is a contradiction that $\Phi$ maps a rank one operator into
 a rank zero operator, so the result follows.
 \end{proof}

The next lemma states a different version of Lemma \ref{tvl} for
the case $\dim N \geq 2$ and
 $\dim N_{-}^{\perp}=1$. The proof of this lemma is almost same
 as the proof of Lemma \ref{tvl}.

\begin{lemma} \label{tvl2}
  Let $\mathfrak{L}$ be a
 commutative subspace lattice and $\Phi$ be a rank preserving map
 on $\alg\mathfrak{L}$. Let $N$ be a subspace in $\alg\mathfrak{L}$
 with $\dim N \geq 2$ and $\dim N_{-}^{\perp}=1$. If we consider the following two equalities,
  \begin{align*}
   \Phi(x \otimes f^{\ast})&= u \otimes v^{\ast} &&\cdots \cdots (G)
      \\
   \Phi(y \otimes f^{\ast})&= p \otimes q^{\ast} &&\cdots \cdots (H)
  \end{align*}
  where $ x,y $ are two linearly independent vectors in $N$ and
  $f \in N_-^{\perp}$, then one of the following holds.
  \begin{enumerate}
    \item $ u \sim p , \, v \nsim q $
    \item $  u \nsim q , \, v \sim q $
   \end{enumerate}
\end{lemma}

With these three lemmas in hand, we now can show the following
important result.
\begin{lemma}\label{ml}
Let $\mathfrak{L}$ be a commutative subspace lattice and $\Phi$ be
a rank preserving map on $\alg\mathfrak{L}$. Let $N$ be a nonzero
element in $\mathfrak{L}$ with $N_{-}\neq I$. Then at least one of
the following holds.
\begin{enumerate}
\item There exist a linear map $U$ from $N$ to $\mathcal{H}$ and a
linear map $V$ from $N_{-}^{\perp}$ to $\mathcal{H}$ such that
$\Phi(x \otimes f^*)=U(x) \otimes V(f)^*$ for all $x \in N$ and $f
\in N_-^{\perp}.$
 \item There exist a conjugate linear map $U$ from
$N_{-}^{\perp}$ to $\mathcal{H}$ and a conjugate linear map $V$
from $N$ to $\mathcal{H}$ such that $\Phi(x \otimes f^*)=U(f)
\otimes V(x)^*$ for all $x \in N$ and $f \in N_-^{\perp}.$
\end{enumerate}
\end{lemma}

\begin{proof}
Suppose $(1)$ holds for all $x,y \in N$ and $f,g \in N_-^{\perp}$
in Lemma \ref{fvl}, Lemma \ref{tvl} and Lemma \ref{tvl2}.

\underline{Case 1:} Suppose that $\dim N=1$ and $\dim
N_{-}^{\perp}=1$. Fix $x_1 \in N, f_1 \in
 N_-^{\perp}$. Set $\Phi(x_1 \otimes f_1^*)=u_1 \otimes
 v_1^*$. For any $x \in N, f \in N_-^{\perp}$, there exist
 $\lambda, \gamma \in \mathbb{C}$ such that $x=\lambda x_1, f=\gamma
 f_1$. If we define $U:N \rightarrow \mathcal{H} $ and $V:
 N_-^{\perp} \rightarrow \mathcal{H}$ by
 $U(x)=\lambda u_1$ and $V(f)=\gamma v_1$, then these maps are
 clearly well defined linear maps. However, if we define $U:N \rightarrow \mathcal{H} $
 and  $V:N_-^{\perp} \rightarrow \mathcal{H}$ by $U(x)=\lambda v_1$ and $V(f)=\gamma
 u_1$, these maps are also well defined linear maps. Therefore
 both conclusions $(1)$ and $(2)$ hold in this special case.

\underline{Case 2:} Suppose that $\dim N=1$ and $\dim
N_{-}^{\perp}\geq 2$. Fix $x_1 \in N$. It follows from Lemma
\ref{tvl} (1) that it is possible to to define a map $V$ on
$N_-^{\perp}$ by $\Phi(x_1 \otimes f^*)=u_1 \otimes V(f) $. For
any $f, g \in N_{-}^{\perp}, t\in \mathbb{C}$, since $\Phi$ is
linear,
 \begin{align*}
 \Phi(x_1\otimes (f+tg)^*)&=\Phi(x_1\otimes
 f^*)+\Phi(x_1\otimes(tg)^*)\\
 &=u_1\otimes V(f)^* +\overline{t}u_1\otimes V(g)^*.
 \end{align*}

From the definition of $V$, $\Phi(x_1\otimes (f+tg)^*)=u_1\otimes
V(f+tg)^*$. By comparing the two equalities, we see the map
$V:N_-^{\perp} \rightarrow \mathcal{H}$
 is a linear map on $N_{-}^{\perp}$.

For any $x \in N$, there exists $\lambda \in \mathbb{C}$ such that
$x=\lambda x_1$. Now we can define a linear map $U:N \rightarrow
 \mathcal{H}$ by $U(x)=\lambda u_1$. Therefore the conclusion $(1)$ holds with this assumption.
Likewise, the conclusion $(2)$ follows if we assume $(2)$ in Lemma
\ref{tvl} is true.

\underline{Case 3:} Suppose that $\dim N \geq 2$ and $\dim
N_{-}^{\perp}= 1$. The proof of this case is almost identical to
the proof of case 2.

\underline{Case 4:}
 Let's fix $x_1\in N$. By Lemma \ref{fvl} (1), we can define a
 map $V$ from $N_{-}^{\perp}$ to $\mathcal{H}$ by $\Phi(x_1 \otimes f^*)=u_1 \otimes
 V(f)^*$. For any $f, g \in N_{-}^{\perp}, t\in \mathbb{C}$, since
 $\Phi$ is linear,
 \begin{align*}
 \Phi(x_1\otimes (f+tg)^*)&=\Phi(x_1\otimes
 f^*)+\Phi(x_1\otimes(tg)^*) \\
 &=u_1\otimes V(f)^* +\overline{t}u_1\otimes V(g)^*.
 \end{align*}

 From the definition of $V$, $\Phi(x_1\otimes (f+tg)^*)=u_1\otimes
 V(f+tg)^*$.
 By comparing the two equalities, we see the map $V$ is a linear
 map on $N_{-}^{\perp}$. Likewise, we can define a linear map $U$
 from $N$ to $\mathcal{H}$ by $\Phi(x\otimes {f_1}^*)=U(x)\otimes
 {v_1}^*$ for some fixed $f_1, v_1 \in N_{-}^{\perp}$. By
 considering
 \begin{align*}
 \Phi(x_1\otimes {f_1}^*)&=U(x_1)\otimes{v_1}^*\\
 &=u_1\otimes V(f_1)^*\\
 &=u_1\otimes{v_1}^* ,
 \end{align*}
 we can make an observation that $U(x_1)=u_1$ and $V(f_1)=v_1$.
 Now for any $x\in N, f\in N_{-}^{\perp}$, we may write
 $\Phi(x\otimes f^*)$ as follows.
 \[ \Phi(x\otimes f^*)=\alpha(x,f)U(x)\otimes V(f)^*\]
 where $\alpha(x,f)$ is a complex valued function.

 We claim that $\alpha(x,f)=1$ for all $x$ and $f$. Note that
 $\Phi(x_1\otimes f^*)=\alpha(x_1,f)u_1\otimes
 V(f)^*=u_1\otimes V(f)^*$, so $\alpha(x_1,f)=1$ for any $f\in
 N_{-}^{\perp}$. Similarly, we observe that $\alpha(x,f_1)=1$ for any
 $x\in N$. Consider
  \begin{align*}
 \Phi((x+sx')\otimes f^*)&=\alpha(x+sx',f)(U(x+sx')\otimes
 V(f)^*)\\
 &=\alpha(x+sx',f)[(U(x)+sU(x'))\otimes V(f)^*]\\
 &=\alpha(x+sx',f)U(x)\otimes V(f)^*+ s\alpha(x+sx',f)U(x')\otimes V(f)^*
  \end{align*}
  where $x$ and $x'$ are linearly independent vectors in $N$.

On the other hand,
  \begin{align*}
  \phi((x+sx')\otimes f^*)&=\Phi(x\otimes f^*)+s\Phi(x'\otimes f^*)\\
  &=\alpha(x,f)U(x)\otimes V(f)^* +s\alpha(x',f)U(x')\otimes
  V(f)^*.
  \end{align*}
 We get the following equality by comparing the two equalities
 above:
\[\alpha(x+sx',f)U(x)+s\alpha(x+sx',f)U(x')=\alpha(x,f)U(x)+s\alpha(x',f)U(x')\]
for all $s\in \mathbb{C}$. Since $U(x)$ and $U(x')$ are linearly
independent, $\alpha(x,f)=\alpha(x+sx',f)=\alpha(x',f)$. Hence,
$\alpha$ is independent of $x$. Likewise we can show that $\alpha$
is independent of $f$. Therefore $\alpha \equiv 1.$ The result
follows.
\end{proof}

What Lemma~\ref{ml} says is that whether a given subspace $N$
satisfies $(1)$ or $(2)$ in the lemma is a characteristic of the
subspace which is associated with the rank preserving map $\Phi$.

\begin{definition}
 Let $\mathfrak{L}$ be a
commutative subspace lattice and $\Phi$ be a rank preserving map
on $\alg\mathfrak{L}$. Let $N$ be a nonzero subspace in $\alg
\mathfrak{L}$. The subspace $N$ is called \emph{consistent with
respect to} $\Phi$ if $\Phi$ and $N$ satisfy $(1)$ in Lemma
\ref{ml}. The subspace $N$ is called \emph{twisted with respect
to} $\Phi$ if $\Phi$ and $N$ satisfy $(2)$ in Lemma \ref{ml}.
\end{definition}

In a nest algebra, it is impossible that both $N$ and $N^{\perp}$
are in a nest $\mathcal{N}$ (a totally ordered lattice) if $N \neq
0$ and $N \neq I.$ But a commutative subspace lattice
$\mathfrak{L}$ can have $N$ and $N^{\perp}$ both in it with the
assumption $N \neq 0$ and $N \neq I.$ If there exists a subspace
$N \in \mathfrak{L}$ with $N_- \neq I$ such that $\dim N=1$ and
$N^{\perp} \in \mathfrak{L}$, then $N_-=\vee\{E \in
\mathfrak{L}\,|\, E \ngeq N\} \geq N^{\perp}$ since $N^{\perp} \in
\mathfrak{L}$ and $N^{\perp} \ngeq N$. Hence $N_-^{\perp} \leq N$.
Since $N_- \neq I$, $N_-^{\perp}=N$. \underline{Case 1} in the
proof of Lemma~\ref{ml} shows that such $N$ is both consistent and
twisted with respect to $\Phi$.

\begin{definition}
A nonzero subspace $N$ in $\mathfrak{L}$ with $\dim N=1$ and
$N^{\perp} \in \mathfrak{L}$ is called \emph{isolated}.
\end{definition}

\begin{lemma} \label{isl}
Let $\Phi$ be a rank preserving map on $\alg\mathfrak{L}$ and $N$
be a subspace in $\mathfrak{L}$. Then $N$ is isolated if and only
if $N$ is both consistent and twisted with respect to $\Phi$.
\end{lemma}
\begin{proof}
It suffices to prove sufficiency. Suppose $N$ is both consistent
and twisted with respect to $\Phi$. Then there exist linear maps
$U_1,V_2$ defined on $N$, and $U_2,V_1$ defined on $N_-^{\perp}$
so that \[\Phi(x \otimes f^*)=U_1x \otimes (V_1f)^*=U_2f \otimes
(V_2x)^*\] \[\Phi(y \otimes f^*)=U_1y \otimes (V_1f)^*=U_2f
\otimes (V_2y)^*\] for arbitrary $x,y \in N$ and $f \in
N_-^{\perp}$. It can be observed that $U_1x \sim U_2f$ and $U_2f
\sim U_1y$, so $U_1x \sim U_1y$. Since the map $U_1$ can not have
nonzero kernel, $x \sim y$. Using the same idea as this with $x
\in N$ and $f,g \in N_-^{\perp}$, we can argue that $f \sim g$.
These facts imply that $\dim N=\dim N_-^{\perp}=1$. Observe that
if $x \in N$, $x \otimes x^*$ is the projection on $N$ since $\dim
N=1$, so $x \otimes x^* \in \alg\mathfrak{L}$. Therefore, $x \in
N_-^{\perp}$. Then $N=N_-^{\perp}$. Thus $N$ is isolated.
\end{proof}

\begin{lemma}\label{sl}
Let $M$ and $N$ be two subspaces in $\mathfrak{L}$ with $N \leq
M$. Let $\Phi$ be a rank preserving map on $\alg\mathfrak{L}$.
Then the followings are true.
\begin{enumerate}
\item If $M$ is consistent with respect to $\Phi$, then $N$ is
also consistent with respect to $\Phi$. \item If $M$ is twisted
with respect to $\Phi$, then $N$ is also twisted with respect to
$\Phi$.
\end{enumerate}
 \end{lemma}

\begin{proof}
First, let $M$ be a consistent element in $\mathfrak{L}$ with
respect to $\Phi$. Without loss of generality, we assume $\dim N
\geq 2$. Suppose that $N$ is twisted with respect to $\Phi$. Let
$U_M,V_M$ and $U_N,V_N$ denote the two maps in lemma \ref{ml}
defined on $M$ and $N$ ,respectively. By the assumption, we can
choose two linearly independent vectors $x,y \in N$ and a vector
$f \in M_-^{\perp}$. Note that $N_-^{\perp} \geq M_-^{\perp}$
since $N \leq M$. Since $M$ is consistent with respect to $\Phi$
and $N$ is twisted with respect to $\Phi$, we have
\[\Phi(x \otimes f^*)= U_Mx \otimes (V_Mf)^*=U_Nf \otimes (V_Nx)^*
\]
   \[ \Phi(y \otimes f^*)= U_My \otimes (V_Mf)^*=U_Nf \otimes (V_Ny)^*. \]
Then there exist two complex numbers $\lambda$ and $\gamma$ such
that $U_Mx=\lambda U_Nf$ and $U_My=\gamma U_Nf$. But this
contradicts the fact that $U_Mx$ and $U_My$ are linearly
independent by Lemma \ref{fvl}, so the result $(1)$ follows. The
proof of $(2)$ is essentially the same as the proof of $(1)$.
\end{proof}

\begin{lemma}\label{vil}
Let $\Phi$ be a rank preserving map on $\alg \mathfrak{L}$. Let
$M$ and $N$ be two non-zero, non-isolated subspaces in
$\mathfrak{L}$ with $M_- \neq I$. If $M$ is consistent with
respect to $\Phi$, and $N$ is twisted with respect to $\Phi$, then
\[M \wedge N=0\] and
\[M_-^{\perp} \wedge N_-^{\perp}=0.\]
\end{lemma}

\begin{proof}
Suppose $M \wedge N \neq 0$. Then $M \wedge N$ is a non-zero
subspace of $M$ and $N$. By Lemma \ref{isl}, $M \wedge N$ is
isolated, so $M \wedge N=(M \wedge N)_-^{\perp}=\langle e \rangle$
for some $e \in M \wedge N$. By Lemma \ref{ul}, we have $(M \wedge
N)_- \leq M_- \wedge N_-$. It follows that $(M \wedge N)_-^{\perp}
\geq M_-^{\perp} \vee N_-^{\perp}$. Thus
$M_-^{\perp}=N_-^{\perp}=(M \wedge N)_-^{\perp}=\langle e
\rangle$. On the other hand, we have $M_- \geq N$ since $N \ngeq
M$, so $\langle e\rangle=M_-^{\perp} \leq N^{\perp}$. Therefore $N
\wedge N^{\perp} \geq \langle e\rangle \neq 0$ which gives a
contradiction. Thus \[M \wedge N=0.\]

 Observe that $M \vee N$ is
neither consistent nor twisted with respect to $\Phi$ by Lemma
\ref{sl}. This implies $(M \vee N)_-=I$ by Lemma \ref{ml}.
Therefore \[M_-^{\perp} \wedge N_-^{\perp}=(M_- \vee
N_-)^{\perp}=(M \vee N)_-^{\perp}=0.\]
\end{proof}

\begin{theorem}\label{dct}
Let $\mathfrak{L}$ be a completely distributive commutative
subspace lattice and $\Phi$ be a rank preserving map on $\alg
\mathfrak{L}$. Let \[M_i=\vee\{E \in \mathfrak{L}\,|\,E \,
\textrm{is isolated}\},\]
    \[M_c=\vee\{E \in \mathfrak{L}\,|\, E \, \textrm{is consistent with
respect to} \, \Phi \} \ominus M_i\] and
    \[M_t=\vee\{E \in \mathfrak{L}\,|\, E \, \textrm{is twisted with respect
    to}\, \Phi \ \} \ominus M_i.\]

Then \[\alg\mathfrak{L}=\alg(M_i\mathfrak{L}M_i) \oplus
\alg(M_c\mathfrak{L}M_c) \oplus \alg(M_t\mathfrak{L}M_t) .\]

\end{theorem}

\begin{proof}
By Lemma \ref{vil}, it is clear that $M_i \oplus M_c \oplus
M_t=I$. Since $M_i^{\perp} \in \mathfrak{L}$, $M_c, M_t \in
\mathfrak{L}.$ Note that $M_c^{\perp}=M_i \vee M_t.$ Hence
$M_c^{\perp} \in \mathfrak{L}$. Likewise, $M_t^{\perp} \in
\mathfrak{L}$. The result follows.
\end{proof}

Remember that the operators $U$ and $V$ were constructed with
 the subspace $N$ fixed. We now refer to
them as $U_N$ and $V_N$ and try to fit together the $U_N$'s and
$V_N$'s into operators $U$ and $V$ defined on the whole space
$\mathcal{H}$ such that
\[\Phi(x \otimes f^*)=Ux \otimes (Vf)^* \]
whenever $x \otimes f^* \in \alg\mathfrak{L}$. Note that if $M_-
 \neq I$, if $M<N$, and if $x \in M$ and $f \in N_-^{\perp}$, then $x \in N$
and $\Phi(x \otimes f^*)=(U_N x) \otimes (V_N f)^*$. On the other
hand, $f \in N_-^{\perp} <M_-^{\perp}$, so $\Phi(x \otimes
f^*)=(U_M x) \otimes (V_M f)^*$. Thus there exists a complex
number $\lambda$ such that $U_M x=\lambda U_N x$ and $V_N
f=\overline{\lambda}V_M f$. Since $x$ and $f$ may vary
independently, $\lambda$ does not depend on $x$ and $f$ but only
on $M$ and $N$. We call it $\lambda_{MN}$ and have
\[ U_M=\lambda_{MN} U_N|M \,\,\, \textrm{and} \,\,\,
V_N=\overline{\lambda_{MN}}V_M|N_-^{\perp} .\]

Let $\mathfrak{F}$ be the collection of all subspaces $N$ in
$\mathfrak{L}$ such that $N \neq 0$ and $N_- \neq I$.

\begin{remark}
In the remainder of this paper we assume that every subspace in
$\mathfrak{F}$ is consistent with respect to $\Phi$.
\end{remark}
 Suppose that $M$ and $N$ lie
in $\mathfrak{F}$. We will say that $M$ and $N$ are
\emph{comparable} if either $M \leq N$ or $N \leq M$. Suppose $M$
and $N$ are comparable. If $M \leq N$, $\lambda_{MN}$ has already
been defined. If $N \leq M$, define
\[ \lambda_{MN}= \frac{1}{\lambda_{NM}}\,\, .\]
Thus $\lambda_{MN}$ is defined whenever $M$ and $N$ are
comparable, and it is easy to check that
$\lambda_{LN}=\lambda_{LM}\lambda_{MN}$ whenever each pair from
$\{L,M,N\}$ is comparable.

\begin{definition}
Let $M$ and $N$ be two subspaces in $\mathfrak{F}$. We define a
\emph{chain} from $M$ to $N$ to be a finite sequence of subspaces
$\{M_0,M_1,...,M_n\}$, each $M_k$ in $\mathfrak{F}$, such that
$M_0=M, M_n=N$, and such that $M_k$ is comparable to $M_{k+1}$ for
each $k=0,1,...,n-1$.

 If $M=N$, the chain $\{M_0,M_1,...,M_n\}$ is
called a \emph{cycle of length} n.
\end{definition}

Suppose $M,N \in \mathfrak{F}.$ If there is a chain
$\{M_0,M_1,...,M_n\}$ from $M$ to $N$, we want to define
$\lambda_{MN}$ to be $\lambda_{MM_1}\lambda_{M_1M_2} \cdots
\lambda_{M_{n-2}M_{n-1}}\lambda_{M_{n-1}N}$. Since there may be
more than one chain from $M$ to $N$, we need to show that such a
product is well defined. The following lemma will prove this. An
elaborate proof of the lemma can be found in \cite{MR87k:47103}.

\begin{lemma} \label{lprt}
Let $\{M_0,M_1,...,M_n\}$ be a cycle in $\mathfrak{F}$. Then
\[ \lambda_{M_0M_1}\lambda_{M_1M_2} \cdots
\lambda_{M_{n-2}M_{n-1}}\lambda_{M_{n-1}M_n}=1  .\]

\end{lemma}

By making use of the Lemma \ref{lprt}, we can now show that
\[\lambda_{MN}=\lambda_{MM_1}\lambda_{M_1M_2} \cdots
\lambda_{M_{n-2}M_{n-1}}\lambda_{M_{n-1}N}\] is well defined if
there is a chain $\{M_0,M_1,...,M_n\}$ from $M$ to $N$.

Suppose that there are two chains
$\{M_0,M_1,M_2,...,M_{n-2},M_{n-1},M_n\}$ and \\
$\{M_0,M_1',M_2',...,M_{n-2}',M_{n-1}',M_n\}$ from $M$ to $N$. We
now form a cycle \[\{M,M_1,...,M_{n-1},N,M_{n-1}',...,M_1',M\}.\]
By applying Lemma \ref{lprt} to this cycle, it follows that
\[\lambda_{MM_1} \cdots \lambda_{M_{n-1}N}\lambda_{NM_{n-1}'}
\cdots \lambda_{M_1'M}=1.\] From this equation we get,
\[\lambda_{MM_1} \cdots \lambda_{M_{n-1}N}=
\frac{1}{\lambda_{NM_{n-1}'}\cdots
\lambda_{M_1'M}}=\lambda_{MM_1'}\cdots \lambda_{M_{n-1}'N}. \]
This shows that $\lambda_{MN}$ is well defined.\\

Recall that $\mathfrak{F}$ denotes the collection of all $N$ in
$\mathfrak{L}$ such that $N \neq 0$ and $N_- \neq I$. Fix $N \in
\mathfrak{L}$ and let
\[\mathfrak{G}_N^n=\{M \in \mathfrak{F}\,|\,M \, \textrm{can be
connected to} \, N \, \textrm{by a chain of length} \, k \leq
n\}.\] Let $\mathfrak{G}_N=\cup_{n}\mathfrak{G}_N^n$.

\begin{lemma}\label{du}
Let $N \in \mathfrak{F}$. Then there exist linear transformations
$U$ with dense domain in $\vee \{M \, | \, M \in \mathfrak{G}_N\}$
and $V$ with dense domain in $\vee \{M_-^{\perp} \, | \, M \in
\mathfrak{G}_N \}$ such that $\Phi(x \otimes f^*)=Ux \otimes
(Vf)^*$ whenever there is $M \in \mathfrak{G}_N$ for which $x \in
M$ and $f \in M_-^{\perp}$.
\end{lemma}

\begin{proof}
For $M \in \mathfrak{G}_N$, we have associated operators $U_M$ and
$V_M$ such that \[\Phi(x \otimes f^*)=U_Mx \otimes (V_Mf)^*,\]
whenever $x \in M$ and $f \in M_-^{\perp}$. Let
$\widetilde{U}_M=\lambda_{NM}U_M$ and
$\widetilde{V}_M=\overline{\lambda}_{MN} V_M$. Since
$\lambda_{MN}\lambda_{NM}=1$, we have $\Phi(x \otimes f^*)=
\widetilde{U}_Mx \otimes (\widetilde{V}_Mf)^*$ for $x \in M$ and
$f \in M_-^{\perp}$. Note that for $L,M \in \mathfrak{G}_N$ with
$L \leq M$, if $\{N,N_1,...,N_{n-1},M \}$ is a chain from $N$ to
$M$, then $\{N,N_1,...,N_{n-1},M,L\}$ is a chain from $N$ to $L$.
Therefore $\lambda_{NL}=\lambda_{NM}\lambda_{ML}$. Thus, if $x \in
L$, we have
$\widetilde{U}_{L}x=\lambda_{NL}U_{L}x=\lambda_{NM}\lambda_{ML}U_{L}x$.
On the other hand, by definition of $\lambda_{LM}$ we have
$U_Lx=\lambda_{LM}U_Mx$, so
$\widetilde{U}_{L}x=\lambda_{NM}\lambda_{ML}\lambda_{LM}U_{M}x=\lambda_{NM}U_{M}x
=\widetilde{U}_{M}x$. Thus, if $L \leq M$, $\widetilde{U}_L$ and
$\widetilde{U}_M$ agree on $L$. Let $\mathcal{M} =
\{x_1+\cdots+x_n \, | \, \textrm{for some positive integer} \, n
\, , \, x_i \in M \, \textrm{for some} \, M \in \mathfrak{G}_N
\}$. Define a linear transformation $U$ on $\mathcal{M}$ by
$U|M=\widetilde{U}_M$. By the coherence of the $\widetilde{U}_M$,
the map $U$ is well defined. Similarly, we can define
$V|M_-^{\perp}=\widetilde{V}_M$. If $x \in M, \, f \in
M_-^{\perp}$ for $M \in \mathfrak{G}_N$, then we have $\Phi(x
\otimes f^*)=Ux \otimes (Vf)^*$.
\end{proof}

\begin{lemma}\label{orth}
For each $N \in \mathfrak{F}$, let $G_N=\vee\{M \, | \, M \in
\mathfrak{G}_N\}$ and $F_N=\vee\{M_-^{\perp} \, | \, M \in
\mathfrak{G}_N \}$. Let $N,N' \in \mathfrak{F}$.
\begin{enumerate}
\item If $N' \in \mathfrak{G}_N$, then $G_N=G_{N'}$ and
$F_N=F_{N'}$
 \item If $N' \notin \mathfrak{G}_N$, then $G_N \perp
G_{N'}$ and $F_N \perp F_{N'}$
\end{enumerate}
\end{lemma}

\begin{proof}
Suppose that $N' \in \mathfrak{G}_N$. Note that if there exists a
chain connecting $N$ to $N'$ and one connecting $N'$ to $M$, then
there is a chain connecting $N$ to $M$, whence
$\mathfrak{G}_N=\mathfrak{G}_{N'}$. Thus $G_N=G_{N'}$ and
$F_N=F_{N'}$. Suppose that $N' \notin \mathfrak{G}_N$; then, for
any $M' \in \mathfrak{G}_{N'}$, $M' \notin \mathfrak{G}_N$. For
such an $M'$, consider the projection $NM'$. If $NM' \neq 0$, then
$NM' \in \mathfrak{G}_N \cap \mathfrak{G}_{M'}=\mathfrak{G}_{M'}
\cap \mathfrak{G}_N$. Since this is impossible, it must be that
$NM'=0$. Likewise, if $M \in \mathfrak{G}_N$, $MM'=0$ and hence
$G_NG_{N'}=0$.
\end{proof}

\begin{theorem} \label{mpn}
Let $\mathfrak{L}$ be a completely distributive commutative
subspace lattice and let $\Phi$ be a weakly continuous rank
preserving map defined on $\alg\mathfrak{L}$. If we assume that
$\mathfrak{F}=\mathfrak{G}_N$ for some $N \in \mathfrak{F}$ (i.e.
$\alg \mathfrak{L}$ is irreducible.), then there exist two densely
defined linear transformations $U,V$ such that
\[\Phi(A)=UAV^*\] for all $A \in \alg\mathfrak{L}$.
\end{theorem}

\begin{proof}
By Lemma \ref{du}, there exist densely defined maps $U$ and $V$
such that $\Phi(x \otimes f^*)=Ux \otimes (Vf)^*$ whenever $x
\otimes f^* \in \alg\mathfrak{L}$. Since
$\mathfrak{F}=\mathfrak{G}_N$ for some $N \in \mathfrak{F}$, the
domain of $U$ is the nonclosed linear span of $\mathfrak{F}=\{E
\in \mathfrak{L}\,|\, E_-\neq I \}$ and the domain of $V$ is the
nonclosed linear span of $\{M_-^{\perp} \,|\, M \in \mathfrak{F}
\}$. By Lemma \ref{cdl1}, these domains are dense in
$\mathcal{H}$. Recall that $R_{\mathfrak{L}}$ denote the linear
span of the rank one operators in $\alg\mathfrak{L}$. Therefore,
for any operator $A \in R_{\mathfrak{L}}$, $\Phi(A)=UAV^*$ by the
linearity of $\Phi$. By Lemma \ref{wk}, in CDC algebra it is
guaranteed that there exists a net of operators in
$R_{\mathfrak{L}}$ weakly converging to any operator in
$\alg\mathfrak{L}$. Let $A$ be an operator in $\alg\mathfrak{L}$.
Then there exists a net $\{A_{\alpha}\}_{\alpha \in J}$ such that
$A_{\alpha}$ converges weakly to $A$. For each $h \in
\operatorname{dom}V^*$ and $k \in \operatorname{dom}U^*$,
\begin{align*}
\langle UA_{\alpha}V^*h,k \rangle &= \langle A_{\alpha}V^*h, U^*k
\rangle  \\    &\rightarrow \langle AV^*h,U^*k \rangle \,\,
\textrm{since} \, A_{\alpha} \, \textrm{weakly converges to} \, A
\\  &= \langle UAV^*h,k \rangle .
\end{align*}
Thus \[\Phi(A_{\alpha})=UA_{\alpha}V^*
\overset{w}{\longrightarrow} UAV^*.
\]
On the other hand, $\Phi(A_{\alpha})$ converges weakly to
$\Phi(A)$ since $\Phi$ is weakly continuous. Thus $\Phi(A)=UAV^*$
for $A \in \alg\mathfrak{L}$.

\end{proof}

Now we are finally in a position to state and prove the main
result of this paper.
\begin{theorem}
Let $\mathfrak{L}$ be a completely distributive commutative
subspace lattice and let $\Phi$ be a weakly continuous rank
preserving map defined on $\alg\mathfrak{L}$. If we assume that
every subspace is consistent with respect to $\Phi$, then there
exist two densely defined linear transformations $U,V$ such that
\[\Phi(A)=UAV^*\] for all $A \in \alg\mathfrak{L}$.
\end{theorem}
\begin{proof}
Let $N \in \mathfrak{F}$ and consider the collection
$\mathfrak{G}_N$. If $\mathfrak{G}_N = \mathfrak{F}$, then this is
nothing but the case of Proposition \ref{mpn}. Suppose there
exists $M \in \mathfrak{F}$ with $M \notin \mathfrak{G}_N$ such
that $G_N G_M=G_M G_N=0$. In this way we can form a sequence
$\{G_{N_i}\}$ of mutually orthogonal projections in
$\mathfrak{L}$. We shall suppress the $N$ and write simply $G_i$.
The separability of $\mathcal{H}$ guarantees that there are no
more than countably many $G_i$'s. We have $\vee G_i=I$ since $\vee
\{N \,:\,N \in \mathfrak{F}\}=I$, whence $G_j^{\perp}=\vee_{i \neq
j}G_i$ is also in $\mathfrak{L}$. Therefore, the algebra
$\alg\mathfrak{L}$ can be written as the direct sum $\sum_i \oplus
\alg(G_i\mathfrak{L}G_i)$. For any $A \in \alg\mathfrak{L}$, we
can write $A=\sum_i \oplus A_i$ where $A_i \in
\alg(G_i\mathfrak{L}G_i).$ By Proposition~\ref{mpn}, there exist
two densely defined linear maps $U_i$ and $V_i$ such that
$\Phi(A_i)=U_iA_iV_i^*$ for each $i$. If we define $U=\sum_i
\oplus U_i$ and $V=\sum_i \oplus V_i$, then $U$ and $V$ are
densely defined on $\mathcal{H}$. Then the result follows.
\end{proof}

\begin{corollary}
Let $\mathcal{H}$ be a finite dimensional Hilbert space and let
$\mathfrak{L}$ be a commutative subspace lattice consisting of
subspaces in $\mathcal{H}.$ Let $\Phi$ be a rank preserving map on
$\alg \mathfrak{L}$. If we assume that every element in
$\mathfrak{L}$ is consistent with respect to $\Phi$, then there
exists a map $\Psi$ defined by $\Psi(A)=A\Phi^{-1}(I)$ such that
$\Psi \circ \Phi$ is an isomorphism on $\alg \mathfrak{L}.$
\end{corollary}

\begin{proof}
By Main Theorem, there are two linear maps $U$ and $V$ such that
$\Phi(A)=UAV^*$. Since $\mathcal{H}$ is finite dimensional, it is
true that $\Phi(I)=UV^*$ is invertible. Therefore
$\Psi(A)=A\Phi^{-1}(I)=A(UV^*)^{-1}$ is a well defined map. Then
\begin{align*}
(\Psi \circ \Phi)(A)&=\Phi(A)(V^*)^{-1}U^{-1}
\end{align*}
\begin{align*}
                    &=UAV^*(V^*)^{-1}U^{-1}  \\
                    &=UAU^{-1}
\end{align*}
Thus $\Psi \circ \Phi$ is an isomorphism implemented by $U.$
\end{proof}

\section{EXAMPLES}
In this section we present some examples of rank preserving maps.
Most of the examples we discuss in this section will be either
$\mathcal{A}_{2n}$ or $\mathcal{A}_{\infty}$. The precise
definition of these algebras is given in \cite{MR85g:47062}. But
for our discussion it is enough to say that the algebras
$\mathcal{A}_{2n}$ are tridiagonal matrices, of size $2n \times
2n$, of the form
\[
\begin{bmatrix}
* & * &   &   &   &   &   &   &   & * \\
  & * &   &   &   &   &   &   &   & \\
  & * & * & * &   &   &   &   &   & \\
  &   &   & * &   &   &   &   &   & \\
  &   &   & * & * & * &   &   &   & \\
  &   &   &   &   & \ddots & &   &   & \\
  &   &   &   &   &   &   & * &   & \\
  &   &   &   &   &   &   & * &   & \\
  &   &   &   &   &   &   & * & * & * \\
  &   &   &   &   &   &   &   &   & *
\end{bmatrix}
\]
where all nonstarred entries are $0$. It can be observed that the
associated lattice consists of certain diagonal projections, hence
it is commutative and completely distributive. Thus the algebra
$\mathcal{A}_{2n}$ is reflexive. The algebra
$\mathcal{A}_{\infty}$ consists of infinite matrices of the form
\[
\begin{bmatrix}
* & * &   &   & \\
  & * &   &   & \\
  & * & * & * & \\
  &   &   & * & \\
  &   &   & * & \\
  &   &   &   & \ddots
\end{bmatrix}
.\]

Once again the associated lattice is commutative and completely
distributive, so the algebra $\mathcal{A}_{\infty}$ is reflexive.
\\
\begin{example}
Consider $\mathcal{A}_4$. Define a map $\Phi_1:\mathcal{A}_4
\rightarrow \mathcal{A}_4$ by
\[
 \left[
\begin{matrix}
a & b & 0 & h  \\
0 & c & 0 & 0  \\
0 & d & e & f  \\
0 & 0 & 0 & g
\end{matrix}
\right] \mapsto \left[
\begin{matrix}
e & f & 0 & d  \\
0 & g & 0 & 0  \\
0 & h & a & b  \\
0 & 0 & 0 & c
\end{matrix}
\right] .\] Then it is easy to check that $\Phi_1$ is rank
preserving. Moreover, the map $\Phi_1$ is implemented by
\[
U =
\begin{bmatrix}
0 & 0 & 1 & 0  \\
0 & 0 & 0 & 1  \\
1 & 0 & 0 & 0  \\
0 & 1 & 0 & 0
\end{bmatrix}
.\] In other words, the map $\Phi_1$ is an isomorphism such that
\[\Phi_1(A)=UAU^{-1}.\] If we define a map $\Phi_2:\mathcal{A}_4
\rightarrow \mathcal{A}_4$ by
\[
 \left[
\begin{matrix}
a & b & 0 & h  \\
0 & c & 0 & 0  \\
0 & d & e & f  \\
0 & 0 & 0 & g
\end{matrix}
\right] \mapsto \left[
\begin{matrix}
g & f & 0 & h  \\
0 & e & 0 & 0  \\
0 & d & c & b  \\
0 & 0 & 0 & a
\end{matrix}
\right] ,\] then $\Phi_2$ is also a rank preserving map. But the
map $\Phi_2$ is not an isomorphism. Instead, we can write the map
$\Phi_2$ as
\[\Phi_2(A)=VA^{T}V^{-1}\] where
\[
V =
\begin{bmatrix}
0 & 0 & 0 & 1  \\
0 & 0 & 1 & 0  \\
0 & 1 & 0 & 0  \\
1 & 0 & 0 & 0
\end{bmatrix}
.\]

\end{example}

\begin{example}
Let $\mathcal{H}$ be a separable Hilbert space and
$\mathcal{B}=\{e_k\}_{k=1}^{\infty}$ be an orthonormal basis for
$\mathcal{H}$. Let
$\mathcal{N}=\{0,[e_1],[e_1,e_2],[e_1,e_2,e_3],....\}$. Then
$\mathcal{N}$ is a nest and the corresponding nest algebra is
\[
\alg\mathcal{N} =
\begin{bmatrix}
* & * & * & * & * & \cdots \\
0 & * & * & * & * & \cdots \\
0 & 0 & * & * & * & \cdots \\
0 & 0 & 0 & * & * & \cdots \\
0 & 0 & 0 & 0 & * & \cdots \\
\vdots  & \vdots & \vdots & \vdots & \vdots & \ddots
\end{bmatrix}
.\]

Let $S$ denote the unilateral shift operator such that
$Se_k=e_{k+1}$. If we define $\Phi:\alg\mathcal{N} \longrightarrow
\alg\mathcal{N}$ by $\Phi(A)=SAS^*$, then
\[
\Phi :
\begin{bmatrix}
a_{11} & a_{12} & a_{13} & a_{14} & a_{15} & \cdots \\
0 & a_{22} & a_{23} & a_{24} & a_{25} & \cdots \\
0 & 0 & a_{33} & a_{34} & a_{35} & \cdots \\
0 & 0 & 0 & a_{44} & a_{45} & \cdots \\
0 & 0 & 0 & 0 & a_{55} & \cdots \\
\vdots  & \vdots & \vdots & \vdots & \vdots & \ddots
\end{bmatrix}
\mapsto
\begin{bmatrix}
0 & 0 & 0 & 0 & 0 & \cdots \\
0 & a_{11} & a_{12} & a_{13} & a_{14} & \cdots \\
0 & 0 & a_{22} & a_{23} & a_{24} & \cdots \\
0 & 0 & 0 & a_{33} & a_{34} & \cdots \\
0 & 0 & 0 & 0 & a_{44} & \cdots \\
\vdots  & \vdots & \vdots & \vdots & \vdots & \ddots
\end{bmatrix}
.\] From this fact we can observe that $\Phi$ is a rank preserving
map on $\alg\mathcal{N}$ but it is not an onto map.
\end{example}

\begin{example}
Consider the algebra $\mathcal{A}_{\infty}$ which is described at
the beginning of Section 3. Let $\mathfrak{L}$ be the lattice
associated with $\mathcal{A}_{\infty}$. We can easily argue that
for any $M,N \in \mathfrak{L}$, $\mathfrak{G}_M=\mathfrak{G}_N$.
Therefore $\mathcal{A}_{\infty}$ is irreducible. Let $\Phi$ be a
map on $\mathcal{A}_{\infty}$ defined by
\[
\Phi:
\begin{bmatrix}
a & b & 0 & 0 & 0 & 0 & \cdots\\
0 & c & 0 & 0 & 0 & 0 & \\
0 & d & e & f & 0 & 0 & \\
0 & 0 & 0 & g & 0 & 0 & \\
0 & 0 & 0 & h & i & j & \\
 \vdots & & & &   &   & \ddots
\end{bmatrix}
\longrightarrow
\begin{bmatrix}
a & \frac{1}{2} b & 0 & 0 & 0 & 0 & \cdots\\
0 & c & 0 & 0 & 0 & 0 & \\
0 & \frac{3}{2} d & e & \frac{3}{4} f & 0 & 0 & \\
0 & 0 & 0 & g & 0 & 0 & \\
0 & 0 & 0 & \frac{5}{4} h & i & \frac{5}{6} j & \\
 \vdots & & & &   &   & \ddots
\end{bmatrix}
.\] Then it is obvious that the map $\Phi:\mathcal{A}_{\infty}
\rightarrow \mathcal{A}_{\infty}$ preserves rank and it can also
be observed that the map $\Phi$ is implemented by an unbounded
operator
\[ U=
\begin{bmatrix}
1 & 0 & 0 & 0 & 0 & \cdots\\
0 & 2 & 0 & 0 & 0 & \\
0 & 0 & 3 & 0 & 0 & \\
0 & 0 & 0 & 4 & 0 & \\
0 & 0 & 0 & 0 & 5 & \\
 \vdots & & & &   & \ddots
\end{bmatrix}
.\] In other words, \[\Phi(A)=UAU^{-1}\] for all $A \in \alg
\mathfrak{L}.$

\end{example}

\bibliographystyle{amsalpha}

\end{document}